\documentstyle[12pt]{article}
\newcommand{\be}[1]{\begin{equation} \label{#1} }
\newcommand{\bea}[1]{\begin{eqnarray} \label{#1} }
\newcommand{\bfi}{\begin{figure}}
\newcommand{\efi}{\end{figure}}
\newcommand{\ee}{\end{equation}}
\newcommand{\eea}{\end{eqnarray}}
\newtheorem{theorem}{Theorem}
\newtheorem{remark}{Remark}
\newtheorem{definition}{Definition}
\newtheorem{example}{Example}
\def\bea{\begin{eqnarray}}
\def\eea{\end{eqnarray}}
\def\ba{\begin{array}}
\def\ea{\end{array}}

\def\D{\Delta}

\def\j{\tilde{j}}
\begin{document}
.

\hspace{3.5cm}

\begin{center}
{\bf ON PRINCIPAL MINORS OF BEZOUT MATRIX} \vspace{1cm}
\end{center}

\hspace{0.5cm}

\begin{center}
{R.G. Airapetyan}

Department of Mathematics, Kettering University, Flint, USA.\\ E-mail :
rhayrape@kettering.edu
\end{center}

\vspace{1cm}

\begin{abstract}
Let $x_1,\dots,x_{n}$ be real numbers,
$P(x)=p_n(x-x_1)\cdot\cdot\cdot(x-x_n)$, and $Q(x)$ be a polynomial of
degree less than or equal to $n$. Denote by $\Delta(Q)$ the matrix of
generalized divided differences of $Q(x)$ with nodes $x_1,\dots,x_n$
and by $B(P,Q)$ the Bezout matrix (Bezoutiant) of $P$ and
$Q$. A relationship  between the corresponding principal
minors, counted from the right-hand lower corner, of the matrices $B(P,Q)$ and $\Delta(Q)$
is established. It implies that if the principal
minors of the matrix of divided differences of a function $g(x)$ are positive
or have alternating signs then the roots of the Newton's interpolation polynomial of $g$
are real and separated by the nodes of interpolation.
\end{abstract}



\section{Introduction.}

In this paper a relationship between two well known matrices is
established. The first one is a Bezout matrix $B$
playing an important role in the theory of separation of
polynomial roots. The second one is Newton's matrix of divided
differences $\Delta$ or, in the case of multiple nodes, Hermite's
matrix of generalized divided differences, playing an important
role in numerical analysis and approximation theory.
In this paper we show that the corresponding principal
minors of $B$ and $\Delta$ counted from the right-hand lower corner are
related by a simple formula (are equal when $p_n=1$). An alternative proof of this result 
can be obtained from the results of \cite{dg}. As a simple application of the
relationship between $B$ and $\Delta$, a theorem about locations
of the roots of interpolation polynomials in terms of the
principal minors of $\Delta$ is established. Many applications of Bezout matrix
can be found in  \cite{syl}-\cite{ahks}.
The results of this paper
were announced without proofs in \cite{air2}.

\section{Main Results.}

With the polynomials $P(x)=\sum\limits_{j=0}^np_jx^j$ and
$Q(x)=\sum\limits_{j=0}^nq_jx^j$ let us associate the bilinear
form
\begin{equation} \label{bdef} \sum\limits_{i,j=1}^n
b_{ij}x^{i-1}y^{j-1}=\frac{P(x)Q(y)-P(y)Q(x)}{x-y},
\end{equation}
which Sylvester \cite{syl} named "Bezoutiant". If the degree of $Q$
is less than the degree of $P$, that is
$Q(x)=\sum\limits_{j=0}^mq_jx^j$, $m<n$, then one adds zero
coefficients $q_{m+1},\ldots,q_n$ to $Q$. In what follows we
assume that $m\le n$ and denote $B(P,Q)=||b_{ij}||_{i,j=1,\dots,n}$.

It has been shown (see \cite{l}) that
$$ \label{bf1}B(P,Q)=\pmatrix{
p_1&p_2&\cdot&\cdot&\cdot&p_n\cr
                    p_2&\cdot&\cdot&\cdot&\cdot &0\cr
                    \cdot&\cdot&&&\cdot&\cdot\cr
                    \cdot&\cdot&&\cdot&&\cdot\cr
                    \cdot&\cdot&\cdot&&&\cdot\cr
                     p_n&0&\cdot&\cdot&\cdot&0}
                     \pmatrix{q_0&\cdot&\cdot&\cdot&q_{n-2}&q_{n-1}\cr
                    0&\cdot&\cdot&\cdot&\cdot &q_{n-2}\cr
                    \cdot&\cdot&&&\cdot&\cdot\cr
                    \cdot&&\cdot&&\cdot&\cdot\cr
                    \cdot&&&\cdot &q_0&\cdot\cr
                     0&\cdot&\cdot&\cdot&0&q_0 }$$
 \begin{eqnarray} -\pmatrix{q_1&q_2&\cdot&\cdot&\cdot&q_n\cr
                    q_2&\cdot&\cdot&\cdot&\cdot &0\cr
                    \cdot&\cdot&&&\cdot&\cdot\cr
                    \cdot&\cdot&&\cdot&&\cdot\cr
                    \cdot&\cdot&\cdot&&&\cdot\cr
                     q_n&0&\cdot&\cdot&\cdot&0}
                     \pmatrix{p_0&\cdot&\cdot&\cdot&p_{n-2}&p_{n-1}\cr
                    0&\cdot&\cdot&\cdot&\cdot &p_{n-2}\cr
                    \cdot&\cdot&&&\cdot&\cdot\cr
                    \cdot&&\cdot&&\cdot&\cdot\cr
                    \cdot&&&\cdot &p_0&\cdot\cr
                     0&\cdot&\cdot&\cdot&0&p_0
                     }. \label{bf1}\end{eqnarray}

The main properties of the Bezoutiant
are (see \cite{kn,l,fd,oo})\,:

\begin{itemize}
  \item
    The defect of the Bezoutiant equals the degree of the greatest
    common divisor of the polynomials $P$ and $Q$.
  \item
    The rank of the Bezoutiant matrix equals the degree of the last
    principal minor of the matrix $B=||b_{i,j}||_{i,j=1,\dots,n}$ which
    does not vanish if, in constructing the consecutive major
    minors, one starts from the lower right-hand corner.
  \item
    If the Bezoutiant matrix is positive definite then both polynomials $P(x)$ and $Q(x)$ have real, distinct
    roots.  Moreover, the roots of $P(x)$ and $Q(x)$ interlace.
  \item
    If all consecutive principal minors starting from the lower right-hand
    corner are positive or have alternating signs, then the roots of $P(x)$ and $Q(x)$ are real, distinct, and interlace.
  \end{itemize}

Since principal minors of Bezoutiants play so important a role, it
seems interesting to find explicit formulas for them. If the roots
$x_1,x_2,\dots,x_n$ of $P(x)$ are simple, such formulas were
established in \cite{air1}\,.

\begin{theorem}
Let $|b_{ij}|^n_{i,j=k+1}$ be the principal minors  counted from
the lower right corner of the Bezoutiant $B(P,Q)$ of
polynomials $P(x)=p_n(x-x_1)\dots(x-x_n)$ and $Q(x)$. Then,
\[|b_{i,j}|_{i,j=k+1}^n=p_n^{2(n-k)}\sum\limits_{\begin{array}{ll}(i_1,\dots,i_{n-k})\subset
(1,\dots,n)\\i_1<i_2<\dots<i_{n-k}\end{array}}\frac{Q(x_{i_1})\cdot\cdot\cdot
Q(x_{i_{n-k}})}{P'(x_{i_1})\cdot\cdot\cdot P'(x_{i_{n-k}})}
\]
\begin{equation}\times\prod\limits_{\begin{array}{ll}(j_1,j_2)\subset
(i_1,\dots,i_{n-k})\\j_1<j_2\end{array}}(x_{j_1}-x_{j_2})^2.
\label{t1} \end{equation}
\end{theorem}

\begin{remark}
If $k=n-1$ then formula (\ref{t1}) becomes
\begin{equation}b_{n,n}=p_n^{2}\sum\limits_{i=1}^n\frac{Q(x_{i})}{P'(x_{i})}.
\label{r1} \end{equation}
\end{remark}

\begin{remark}
Since $|b_{ij}|^m_{i,j=k+1}$ are continuous functions of
$x_1,\dots,x_{m}$ in case of multiple roots one has to find the
corresponding limit which is technically difficult and leads to
complicated expressions.
\end{remark}

In order to consider the case of $x_1,\dots,x_{n}$  which are not
necessarily different, let us introduce the following generalized
divided differences.

\begin{definition} (see \cite{sb}):
\[
g[x_i]:=g(x_i),\quad i=1,\dots,n,
\]
\begin{equation} \label{gfd}
g[x_{i_1},\dots,x_{i_k}]:=\left\{\begin{array}{lll}\frac{g[x_{i_2},\dots,x_{i_k}]-g[x_{i_1},
\dots,x_{i_{k-1}}]}{x_{i_k}-x_{i_1}},\,
\hbox{ if }\, x_{i_1}\neq x_{i_k}
\\ \\ \frac{d}{dx}g[x,x_{i_2},\dots,x_{i_{k-1}}]_{x=x_{i_1}}, \,
\hbox{ if }\, x_{i_1}=x_{i_k}.\end{array}\right.
\end{equation}\end{definition}

\begin{remark}
This definition of generalized divided differences is equivalent to the definition given in \cite{sb} if $x_1\le x_2\le\dots\le x_n$.
\end{remark}

Consider the following triangular matrix of the generalized
divided differences: $\Delta(g)=||\Delta_{ij}||_{i,j=1,\dots,n}$, where
\begin{equation} \label{ddm}
\Delta_{ij}=\left\{\begin{array}{ll}0,\,\,\,\hbox{if}\,\,\, i+j<n+1\\
g[x_{n-i+1},\dots,x_j],\,\,\,\hbox{if}\,\,\, i+j\ge n+1,\end{array}\right.
\end{equation}
that is

 \begin{equation} \label{delmat}\Delta(g)=\pmatrix{&&&&&0&\Delta_n\cr
 &0&&&.&\Delta_{n-1}&\Delta_{n-1,n}\cr
 &&.&&.&.&.\cr
 &.&.&&&.&.\cr
 0&.&&&&.&.\cr
 \Delta_1&\Delta_{1,2}&.&.&.&\Delta_{1,n-1}&\Delta_{1,n}}.
 \end{equation}

\begin{remark}
As it is well known, Newton-Hermite's interpolation polynomial for $n$ nodes $\{x_1,\dots,x_n\}$ is
$\Delta_1+\Delta_{1,2}(x-x_1)+\dots+\Delta_{1,n}(x-x_1)\dots(x-x_{n-1})$.
 \end{remark}

 Denote by $|\Delta_{i,j}|_{i,j=k+1}^n$ the principal minors of the matrix $\Delta$
 counted from the lower right corner. The following theorem
 establishes a relationship between principal minors of the
 Bezoutiant and Newton's matrix.

\begin{theorem}\label{t2}
Let $|b_{i,j}|_{i,j=k+1}^n$ and $|\Delta_{i,j}|_{i,j=k+1}^n$ be
the principal minors of the matrices $B(P,Q)$ and $\Delta(Q)$ counted from
the lower right corner. Then
\begin{equation}\label{tf2}
|b_{i,j}|_{i,j=k+1}^n=p_n^{n-k}|\Delta_{i,j}|_{i,j=k+1}^n,\quad
k=0,1,\dots,n-1.\end{equation}\end{theorem}

The relationship between $B(P,Q)$ and $\Delta(Q)$ established in this
theorem is surprising taking into account that these matrices are of
very different type, the Bezoutiant is a symmetric matrix and
Newton's matrix is a triangular matrix. Two simple examples below
show these matrices for some polynomials of degree three.

\begin{example}
Let us consider polynomials $P(x)=x^3-4x^2-x+4=(x+1)(x-1)(x-4)$
and $Q(x)=x^3-6x^2+11x-6=(x-1)(x-2)(x-3)$. Then
\[B(P,Q)=\pmatrix{-38&48&-10\cr
48&-60&12\cr-10&12&-2},\quad
\Delta(Q)=\pmatrix{0&0&6\cr0&0&2\cr-24&12&-2}.\]
Since $p_3=1$ the corresponding principal minors of these two
matrices counted from the lower right-hand corner are equal, they
are $-2,-24,0$.
\end{example}

\begin{example}
Consider polynomials $P(x)=x^3-12x^2+44x-48=(x-2)(x-4)(x-6)$ and
$Q(x)=x^3-9x^2+23x-15=(x-1)(x-3)(x-5)$. Then
\[B(P,Q)=\pmatrix{444&-252&33\cr-252&153&-21\cr33&-21&3},\quad
\Delta(Q)=\pmatrix{0&0&15\cr0&-3&9\cr3&-3&3}.\]
Principal minors counted from the lower right-hand corner are
$3,18,135$.
\end{example}

Theorem \ref{t2} and the properties of the Bezoutiant
described above imply the following theorem.

\begin{theorem} If all consecutive principal minors of the matrix of divided differences (see (\ref{delmat})) of some function $g(x)$ starting from the lower
right-hand corner
 are positive or have alternating signs, then the roots of Newton's interpolation
 polynomial are real, distinct, and interlace with the nodes of interpolation.\end{theorem}

\section{Proofs.}

As it is shown in \cite{l},
$$
(-1)^\frac{(n-k)(n-k-1)}{2}|b_{i,j}|_{i,j=k+1}^n
$$
$$
=\left|\begin{array}{ccccccccccccc}p_n&.&.&.&.&p_{k+1}&|&
p_k&.&.&.&.&p_{2k-n+1}\cr
0&p_n&.&.&.&p_{k+2}&|&p_{k+1}&.&.&.&.&p_{2k-n+2}\cr
.&.&.&&&.&|&.&&.&&&.\cr .&&.&.&&.&|&.&&&.&&.\cr
.&&&.&.&.&|&.&&&&.&.\cr
0&.&.&.&0&p_n&|&p_{n-1}&.&.&.&.&p_k\cr-&-&-&-&-&-&-&-&-&-&-&-&-\cr
q_n&.&.&.&.&q_{k+1}&|&q_k&.&.&.&.&q_{2k-n+1}\cr
0&q_n&.&.&.&q_{k+2}&|&q_{k+1}&.&.&.&.&q_{2k-n+2}\cr
.&.&.&&&.&|&.&&.&&&.\cr .&&.&.&&.&|&.&&&.&&.\cr
.&&&.&.&.&|&.&&&&.&.\cr
0&.&.&.&0&q_n&|&q_{n-1}&.&.&.&.&q_k\end{array}\right|.
$$
Since $Q$ is a polynomial of degree $m$, $ |b_{i,j}|_{i,j=k+1}^n=0$ for $k=m+1,\dots,n-1$. Thus, let
us assume that $k\leq m$ and $k\leq n-1$. Then,
\be{min12}
(-1)^\frac{(n-k)(n-k-1)}{2}|b_{i,j}|_{i,j=k+1}^n=p_n^{n-m}d,
\ee
where
$$
d=\left|\begin{array}{cccccccccccc}p_n&.&.&.&.&.&.&p_{2k-m+1}\cr
0&p_n&.&.&.&.&.&p_{2k-m+2}\cr -&-&-&-&-&-&-&-\cr
-&-&-&-&-&-&-&-\cr 0&.&.&.&0&p_n&\ldots&p_{k}\cr
q_m&.&.&.&.&.&.&q_{2k-n+1}\cr 0&q_m&.&.&.&.&.&q_{2k-n+2}\cr
-&-&-&-&-&-&-&-\cr -&-&-&-&-&-&-&-\cr
0&.&.&.&0&q_{m}&\ldots&q_k\end{array}\right|.
$$
This determinant can be represented as,
 \be{min2}
d=\left|\begin{array}{cccccccccccccc}p_n&.&.&.&p_{n+k-m+1}&.&.&.&p_0&0&.&.&0\cr
.&.&&&.&&&&&.&&&.\cr .&&.&&.&&&&&&.&&.\cr .&&&.&.&&&&&&&.&0\cr
0&.&.&.&p_n&.&.&.&.&.&.&.&p_0\cr-&-&-&-&-&-&-&-&-&-&-&-&-\cr
q_m&.&.&.&.&.&.&.&q_0&0&.&.&0\cr .&.&&&&&&&&.&&&.\cr
.&&.&&&&&&&&.&&.\cr .&&&.&&&&&&&&.&0\cr
0&.&.&.&q_m&.&.&.&.&.&.&.&q_0\cr-&-&-&-&-&-&-&-&-&-&-&-&-\cr
0&.&.&.&.&.&0&|&&&&&&\cr
0&.&.&.&.&.&0&|&&&I_k&&&\cr0&.&.&.&.&.&0&|&&&&&&\end{array}\right|,\ee
where $I_k$ is $k\times k$ unit matrix (obviously, there are no rows below the second dashed line if $k=0$).

First let us assume that the roots of the polynomials $P$ and $Q$ are simple and distinct.

Denote by $V_j(x_1,\dots,x_n,y_1,\dots,y_{n-k})$ the following matrix:
$$\pmatrix{x_1^{j}&\ldots&x_n^{j}&y_1^{j}&\ldots&y_{m-k}^{j}\cr
\vdots&\ldots&\vdots&\vdots&\ldots&\vdots\cr
1&\ldots&1&1&\ldots&1}
$$
Then $V_{n+m-k-1}(x_1,\dots,x_n,y_1,\dots,y_{n-k})$ is the Vandermont matrix and
$$
\det(V_{n+m-k-1}(x_1,\dots,x_n,y_1,\dots,y_{m-k}))
$$
$$=\prod\limits_{1\le i_1<i_2\le n}(x_{i_1}-x_{i_2})
\prod\limits_{\begin{array}{cc}1\le i_1\le n\\1\le i_2\le
m-k\end{array}}(x_{i_1}-y_{i_2})
 \prod\limits_{1\le i_1<i_2\le m-k}(y_{i_1}-y_{i_2})
$$
\be{vand} =\frac{(-1)^{n(m-k)}}{p_n^{m-k}}f(y_1)\ldots
f(y_{m-k})\prod\limits_{1\le i_1<i_2\le n}(x_{i_1}-x_{i_2})
\prod\limits_{1\le j_1<j_2\le m-k}(y_{j_1}-y_{j_2}).\ee
Multiplying determinants(\ref{min2}) and (\ref{vand}) one gets:
$$
d\cdot Vand=\left|\begin{array}{ccc}0&M_1\cr M_2&0\cr
V_{k-1}(x_1,\dots,x_n)&V_{k-1}(y_1,\dots,y_{m}),
\end{array}\right|$$
where
$$
M_1=\pmatrix{y_1^{m-k-1}P(y_1)&\ldots&y_{m-k}^{m-k-1}P(y_{m-k})\cr
\vdots&&\vdots\cr
P(y_1)&\ldots&P(y_{m-k})},
$$
$$
M_2=\pmatrix{x_1^{n-k-1}Q(x_1)&\ldots&x_n^{n-k-1}Q(x_n)\cr
\vdots&&\vdots&\cr
Q(x_1)&\ldots&Q(x_n)}.
$$
Therefore,
\be{bigdet1}
d\cdot Vand=(-1)^{n(m-k)}P(y_1)\dots
P(y_{m-k})D\prod\limits_{1\le j_1<j_2\le m-k}(y_{j_1}-y_{j_2}),\ee
where
 \be{det3} D=\left|\begin{array}{ccccccc}
x_1^{n-k-1}Q(x_1)&\ldots&x_n^{n-k-1}Q(x_n)\cr \vdots&&\vdots\cr
Q(x_1)&\ldots&Q(x_n)\cr x_1^{k-1}&\ldots&x_n^{k-1}\cr
\vdots&&\vdots\cr 1&\ldots&1
\end{array}\right|,\hbox{ if }k\geq 1,
\ee
and
 \be{det31} D=\left|\begin{array}{ccccccc}
x_1^{n-1}Q(x_1)&\ldots&x_n^{n-1}Q(x_n)\cr \vdots&&\vdots\cr
Q(x_1)&\ldots&Q(x_n)
\end{array}\right|,\hbox{ if }k=0.
\ee
From (\ref{min12}), (\ref{vand}), and (\ref{bigdet1}) one gets
\be{det2}
|b_{i,j}|_{i,j=k+1}^n=\frac{(-1)^{(n-k)(n-k-1)/2}p_n^{n-k}}{\prod\limits_{1\le
i_1<i_2\le n}(x_{i_1}-x_{i_2})}D. \ee

Since the case $k=0$ is trivial, let us assume that $k\geq 1$.
Subtracting from all rows of the matrix $D$, except of the $(n-k)$th row and of the
last row, the next row, multiplied by $x_1$, and from the
$(n-k)$th row the last row multiplied by $Q(x_1)$ one gets
\be{det4} D=\left|\begin{array}{cccccccc}
0&x_2^{n-k-2}(x_2-x_1)Q(x_2)&\ldots&x_n^{n-k-2}(x_n-x_1)Q(x_n)\cr
\vdots&&\vdots\cr 0&(x_2-x_1)Q(x_2)&\ldots&(x_n-x_1)Q(x_n)\cr
0&Q(x_2)-Q(x_1)&\ldots&Q(x_n)-Q(x_1)\cr
0&x_2^{k-2}(x_2-x_1)&\ldots&x_n^{k-2}(x_n-x_1)\cr
\vdots&\vdots&\ldots&\vdots\cr 0&(x_2-x_1)&\ldots&(x_n-x_1)\cr
1&1&\ldots&1
\end{array}\right|.
\ee
After pulling out the common multipliers $x_j-x_1$, $j=2,\dots,n$ from columns, one
obtains: \be{det5} D=(-1)^{n+1}\prod\limits_{2\le j\le
n}(x_{j}-x_{1})\left|\begin{array}{ccccccc}
x_2^{n-k-2}Q(x_2)&\ldots&x_n^{n-k-2}Q(x_n)\cr \vdots&&\vdots\cr
Q(x_2)&\ldots&Q(x_n)\cr Q[x_1,x_2]&\ldots&Q[x_1,x_n]\cr
x_2^{k-2}&\ldots&x_n^{k-2}\cr \vdots&\ldots&\vdots\cr 1&\ldots&1
\end{array}\right|.
\ee

Now let us consider two cases: $n-k\ge k-1$ and $n-k<k-1$.
Denote $Q[i]:=Q(x_i)$, $Q[i_1,i_2,\dots,i_k]:=Q[x_{i_1},x_{i_2},\dots,x_{i_k}]$.

{\bf The first case.} If $n-k\ge k-1$, then, continuing this
process, after $k$ steps one gets
$$
D=(-1)^{\frac{k(2n-k+3)}{2}}\prod\limits_{\begin{array}{cc}1\le
j_1\le k\\j_1<j_2\le
n\end{array}}(x_{j_2}-x_{j_1})
$$
\be{det6}
\times\left|\begin{array}{ccccccc}
x_{k+1}^{n-2k-1}Q[k+1]&\ldots&x_n^{n-2k-1}Q[n]\cr
\vdots&&\vdots\cr Q[k+1]&\ldots&Q[n]\cr
Q[k,k+1]&\ldots&Q[k,n]\cr
Q[k-1,k,k+1]&\ldots&Q[k-1,k,n]\cr \vdots&\ldots&\vdots\cr
Q[1,\dots,k,k+1]&\ldots&Q[1,\dots,k,n]
\end{array}\right|.
\ee
Subtracting from the first $n-2k-1$ rows the next row multiplied
by $x_{k+1}$ one obtains:
\be{det7}
D=(-1)^{\frac{k(2n-k+3)}{2}}\det(C)\prod\limits_{\begin{array}{cc}1\le
j_1\le k+1\\j_1<j_2\le n\end{array}}(x_{j_2}-x_{j_1}),
\ee
where the columns of the matrix $C$ are:
$$ C_1=\pmatrix{0\cr\vdots\cr 0\cr Q[k+1]\cr Q[k,k+1]\cr \vdots\cr Q[1,\dots,k,k+1]},
C_i=\pmatrix{x_{k+i}^{n-2k-2}(x_{k+i}-x_{k+1})Q[k+i]\cr \vdots\cr (x_{k+i}-x_{k+1})Q[k+i]\cr Q[k+i]\cr
Q[k,k+i]\cr Q[k-1,k,k+i]\cr\vdots\cr Q[1,\dots,k,k+2]},
$$
$i=2,\dots, n-k$. (The notation $C$ will be used below to denote different matrices.)

After subtracting the first column from the other
columns and pulling out the common factors
$(x_{k+1}-x_k)$, \dots, $(x_n-x_k)$, one gets
$$
D=(-1)^{\frac{k(2n-k+3)}{2}}\prod\limits_{\begin{array}{cc}1\le
j_1\le k+1\\j_1<j_2\le n\end{array}}(x_{j_2}-x_{j_1})\times
$$
\be{det71} \left|\begin{array}{cccccccc}
0&x_{k+2}^{n-2k-2}Q[k+2]&\ldots&x_n^{n-2k-2}Q[n]\cr
\vdots&\vdots&&\vdots\cr 0&Q[k+2]&\ldots&Q[n]\cr
Q[k+1]&Q[k+1,k+2]&\ldots&Q[k+1,n]&\cr
Q[k,k+1]&Q[k,k+2]&\ldots&Q[k,n]\cr
Q[k-1,k,k+1]&Q[k-1,k,k+2]&\ldots&Q[k-1,k,n]\cr
\vdots&\vdots&\ldots&\vdots\cr
Q[1,\dots,k,k+1]&Q[1,\dots,k,k+2]&\ldots&Q[1,\dots,k,n]
\end{array}\right|.\ee
Again, let us subtract from the first $n-2k-2$ rows the next row
multiplied by $x_{k+2}$, then, subtract the second column from the
other columns, and pull out from the columns common factors
$(x_{k+2}-x_{k+1})$, \dots, $(x_n-x_{k+1})$.
Continuing this process, after $n-2k-1$ steps one gets:
$$
D=(-1)^{\frac{k(2n-k+3)}{2}}\det(C)\prod\limits_{\begin{array}{cc}1\le
j_1\le n-k-1\\j_1<j_2\le n\end{array}}(x_{j_2}-x_{j_1}),
$$
where the columns of $C$ are
$$
C_i=\pmatrix{0\cr\vdots\cr 0\cr Q[k+i] \cr Q[k+i-1,k+i]\cr\vdots\cr Q[1,\dots,k+i]},\quad i=1,\dots,n-2k-1,
$$
$$
C_i=\pmatrix{Q[k+i] \cr Q[n-k-1,k+i]\cr\vdots\cr Q[1,\dots,n-k-1,k+i]},\quad i=n-2k,\dots,n-k.
$$
Let us subtract from the last $k$ columns the previous column and factor
out $\prod\limits_{j=n-k+1}^n(x_j-x_{j-1})$. Then let us repeat this
procedure with  the last $k-1$ columns, and so on. Finally,
$$
D=(-1)^{\frac{k(2n-k+3)}{2}}\prod\limits_{1\le j1<j_2\le
m}(x_{j_2}-x_{j_1})\times
$$
\be{det10} \left|\begin{array}{cccccccccc}
0&\ldots&Q[n-k]&\ldots&Q[x_{n-k},\dots,x_{n}]\cr
\vdots&\ldots&\vdots&\ldots&\vdots\cr
0&\ldots&Q[x_{k+2},\dots,x_{n-k}]&\ldots&Q[x_{k+2},\dots,x_n]\cr
Q[k+1]&\ldots&Q[x_{k+1},\dots,x_{n-k}g&\ldots&Q[x_{k+1},\dots,x_n]\cr
Q[x_{k},x_{k+1}]&\ldots&Q[x_{k},\dots,x_{n-k}]&\ldots&Q[x_{k},\dots,x_n]\cr
\vdots&\ldots&\vdots&\ldots&\vdots\cr
Q[x_{1},\dots,x_{k+1}]&\ldots&Q[x_{1},\dots,x_{n-k}]&\ldots&Q[x_{1},\dots,x_n]
\end{array}\right|.\ee
Since, $\prod\limits_{1\le j_1<j_2\le
n}(x_{j_2}-x_{j_1})=(-1)^{n(n-1)/2}\prod\limits_{1\le j1<j_2\le
n}(x_{j_1}-x_{j_2})$, one obtains (\ref{tf2}) from (\ref{ddm}), (\ref{det2}) and
(\ref{det10}).

\vspace{1cm}

{\bf The second case.} If $n-k<k-1$, then, transforming the
determinant similarly to the first case, from (\ref{det5}) one
obtains:
$$
D=(-1)^{\frac{(n+k+4)(n-k-1)}{2}}\det(C)\prod\limits_{\begin{array}{cc}1\le
j_1\le n-k-1\\j_1<j_2\le
n\end{array}}(x_{j_2}-x_{j_1}),
$$
where the columns of $C$ are
$$
C_i=\pmatrix{Q[n-k-1+i]\cr Q[1,n-k-1+i]\cr\vdots\cr Q[1,\dots,n-k-1,n-k-1+i]\cr
x_{n-k-1+i}^{2k-n} \cr \vdots\cr x_{n-k-1+i}\cr 1},\quad i=1,\dots,k+1.
$$
After subtracting from the first $n-k$ rows the last row multiplied by the first element of the row, from
the $n-k+1$th to $k$th rows the next row, multiplied by $x_{n-k}$, and pulling out the common factors $x_{n-k+1}-x_{n-k}$,\dots, $x_{n}-x_{n-k}$
one gets
$$
D=(-1)^{\frac{(n+k+3)(n-k)}{2}}\det(C)\prod\limits_{\begin{array}{cc}1\le
j_1\le n-k\\j_1<j_2\le
n\end{array}}(x_{j_2}-x_{j_1}),
$$
$$
C_i=\pmatrix{Q[n-k,n-k+i]\cr Q[n-k-1,n-k,n-k+i]\cr\vdots\cr Q[1,\dots,n-k,n-k+i]\cr
x_{n-k+i}^{2k-n-1} \cr \vdots\cr x_{n-k+i}\cr 1},\quad i=1,\dots,k.
$$
Continuing this process, after $2k-n$ steps one obtains:
$$
D=(-1)^{\frac{(k-1)(2n-k+4)}{2}}\det(C)\prod\limits_{\begin{array}{cc}1\le
j_1\le k-1\\j_1<j_2\le
n\end{array}}(x_{j_2}-x_{j_1})
$$
$$
C_i=\pmatrix{Q[n-k,\dots,k-1,k-1+i]\cr Q[n-k-1,\dots,k-1,k-1+i]\cr\vdots\cr Q[1,\dots,k-1,k-1+i]\cr
 1},\quad i=1,\dots,n-k+1.
$$
After subtracting the first column from the next columns and pulling out the common factors $x_{k+1}-x_{k}$,\dots, $x_{n}-x_{k}$
one gets
$$
D=(-1)^{\frac{k(2n-k+3)}{2}}\det(C)\prod\limits_{\begin{array}{cc}1\le
j_1\le k\\j_1<j_2\le n\end{array}}(x_{j_2}-x_{j_1}),
$$
$$
C_i=\pmatrix{Q[n-k,\dots,k,k+i]\cr Q[n-k-1,n-k,\dots,k,k+i]\cr\vdots\cr Q[1,\dots,k,k+i]},\quad i=1,\dots,n-k.
$$
Finally, one has to subtract the first column from the next columns and pull out the common factors $x_{k+2}-x_{k+1}$,\dots, $x_{n}-x_{k+1}$, then one has to subtract the second column from the next columns and pull out the common factors $x_{k+3}-x_{k+2}$,\dots, $x_{n}-x_{k+2}$, and so on. Then one gets
$$
D=(-1)^{\frac{k(2n-k+3)}{2}}\prod\limits_{\begin{array}{cc}1\le
j_1<j_2\le n\end{array}}(x_{j_2}-x_{j_1})
$$
\be{seccase3} \times\left|\begin{array}{ccccccc}
\D_{n-k,k+1}(Q)&\D_{n-k,k+2}(Q)&\ldots&\D_{n-k,n}(Q)\cr
\D_{n-k-1,k+1}(Q)&\D_{n-k-1,k+2}(Q)&\ldots&\D_{n-k-1,n}(Q)\cr
\vdots&\vdots&\ldots&\vdots\cr
\D_{1,k+1}(Q)&\D_{1,k+2}(Q)&\ldots&\D_{1,n}(Q)\cr
\end{array}\right|.
\ee

Then, (\ref{det2}) and (\ref{seccase3}) imply (\ref{tf2}).

Thus, Theorem \ref{t2} has proven in the case of simple and distinct roots of polynomials $P$ and $Q$.

To prove the theorem in the general situation, we show first that it remains true if polynomials have common roots.

Let $P(x)=p_n\prod\limits_{j=1}^n(x-x_j)$ and $Q(x)=q_m\prod\limits_{j=1}^r(x-x_j)\prod\limits_{j=1}^{m-r}(x-y_j)$ for some $r$, $0<r\le\min(m,n)$. Then for sufficiently small $\varepsilon$ polynomials $P(x)$ and $Q_\varepsilon(x)=q_m\prod\limits_{j=1}^r(x-x_j-\varepsilon)\prod\limits_{j=1}^{m-r}(x-y_j)$ have distinct roots and therefore (\ref{tf2}) holds. Since both sides in (\ref{tf2}) are continuous functions of $\varepsilon$, the formula remains true when $\varepsilon\to 0$.

Similarly one can prove (\ref{tf2}) in the case of multiple roots. Denote by $r$ the highest multiplicity of the roots of $P(x)$. We will use induction with respect to $r$. If $r=1$ the roots of $P(x)$ are simple. Assume that  (\ref{tf2}) is true for some $r$ and prove it for $r+1$. Assume that there is one root of multiplicity $r+1$. (The case of several roots of multiplicity $r+1$ can be proved similarly.) Let  $P(x)=(x-x_1)^{r+1}P_1(x)$, where $P_1(x)$ is a polynomial of degree $n-r-1$ with the roots distinct from $x_1$. Let $P_\varepsilon(x)=(x-x_1-\varepsilon)(x-x_1)^{r}P_1(x)$  for sufficiently small $\varepsilon$. By assumption,  (\ref{tf2}) is true for $P_\varepsilon$ and $Q$. The left hand side in (\ref{tf2}) is a continuous function of $\varepsilon$. In the right hand side
$$
\lim\limits_{\varepsilon\to 0} g[x_1+\varepsilon,\underbrace{x_1,\dots,x_1}_r,x_2,\dots,x_s]=g[\underbrace{x_1,\dots,x_1}_{r+1},x_2,\dots,x_s].
$$
 This observation completes the proof of Theorem \ref{t2} in the case of multiple roots.

\end{document}